\documentclass[10pt]{article}
\usepackage{amssymb}
\usepackage{amsmath}
\numberwithin{equation}{section}
\usepackage{graphicx}

\newtheorem{rem}{Remark}[section]
\newtheorem{thm}{Theorem}[section]
\newtheorem{lemma}[thm]{Lemma}
\newtheorem{cor}[thm]{Corollary}
\newtheorem{prop}[thm]{Proposition}

\newtheorem{example}{Example}

\newcommand\qed{\hfill\blacksquare\bigskip}

\newcommand\PP{\ensuremath{\mathbb{P}}}
\newcommand\ZZ{\ensuremath{\mathbb{Z}}}
\newcommand\CC{\ensuremath{\mathbb{C}}}

\newcommand\vect[1]{\mbox{\boldmath$#1$}}
\newcommand\Abel{\widetilde{\vect{A}}}
\newcommand\abel{\vect{A}}
\newcommand\tee{\mathcal{T}}

\newcommand\zet[1]{\left\vert {#1} \right\vert}

\newcommand\proof{{\textbf {Proof.}}\ \,}
\begin{document}

\title
{Solution of the generalized periodic discrete Toda equation II;
Theta function solution}

\author{Shinsuke Iwao\\
Graduate School of Mathematical Sciences,\\
The University of Tokyo,\\
3-8-1
Komaba Meguro-ku, Tokyo \ 153-8914, Japan}


\maketitle

\begin{abstract}
We construct the theta function solution to the initial value problem for the
generalized periodic discrete Toda equation.
\end{abstract}

\section{Introduction}\label{sec1}

The aim of the present paper is to obtain an explicit formula
for the solution to the \textit{hungry periodic discrete Toda equation}
(hpdToda) 
(\ref{toda1}--\ref{toda2.5}):
$\forall n,t\in\ZZ$,
\begin{align}
&I_n^{t+M}=I_n^t+V_n^t-V_{n-1}^{t+1},\label{toda1}\\
&V_n^{t+1}=\frac{I_{n+1}^tV_n^t}{I_n^{t+M}},\label{toda2}\\
&I_n^t= I_{n+N}^t,\quad V_n^t=V_{n+N}^t,\label{toda2.5}
\end{align}
where $N$ and $M$ are positive integers.
$t$ is the time variable and $n$ means the position, 
and relation (\ref{toda2.5}) is just 
the periodic boundary condition.
This system is a variant of the 
periodic discrete Toda equation, which
is the $M=1$ case \cite{Tokihiro}.
\bigskip

This article is a continuation of the paper \cite{iwao}.
We will construct a tau function solution for the hungry periodic 
discrete Toda equation (hpdToda).
\bigskip

{\bf Remark}: \label{rema}
To avoid a non-interesting solution 
$I_n^{t+M}=V_n^t$, $V_n^{t+1}=I_{n+1}^t$,
we should assume the extra constraint
\[\textstyle
\prod_{n=1}^{N}{I_n^{t+M}}=
\prod_{n=1}^{N}{I_n^{t}}\,\neq\,
\prod_{n=1}^{N}{V_n^{t+1}}=
\prod_{n=1}^{N}{V_n^{t}},
\]
which is enough to guarantee the existence of a unique solution.
See theorem \ref{thm}.

{\bf Notation}: For a meromorphic function $f$ over a complete curve $C$,
$(f)_0$ (resp.\,$(f)_\infty$) denotes the divisor of zeros (resp.\,poles)
of $f$. Let
$(f):=(f)_0-(f)_\infty$.
$\mathrm{Div}^d(C)$ means the set of divisors
over $C$ of degree $d$ and
$\mathrm{Pic}^d(C)$ means the quotient set defined by
$\mathrm{Pic}^d(C)=\mathrm{Div}^d(C)/(\mbox{linearly equivalent})$.
For an element $\mathcal{D}\in\mathrm{Div}^d(C)$,
$[\mathcal{D}]$ means the image of $\mathcal{D}$ under the natural map 
$\mathrm{Div}^d(C)\to\mathrm{Pic}^d(C)$.
\bigskip

In sections \ref{sec2} and \ref{sec3}, we consider the case $\mathrm{g.c.d.}(N,M)=1$.
We will discuss
the general cases in section \ref{sec4}.

\section{Linearization of hpdToda} \label{sec2}
We summarize the results of \cite{iwao} briefly in this section. The reader
should consult the paper for further details.

\subsection{The spectral curve and the eigenvector mapping}\label{sec2.1}

The hpdToda equation (\ref{toda1}--\ref{toda2.5}) is equivalent to
the following matrix equation:
\begin{equation}\label{matrixform}
L_{t+1}(y)R_{t+M}(y)=R_t(y)L_t(y),
\end{equation}
where $L_t(y)$ and $R_t(y)$ are given by
\[{}
L_t(y)=
\left(\begin{array}{@{\,}cccc@{\,}}
	1 &  & & V_N^t\cdot 1/y \\
	V_1^t & 1 & &  \\
	 & \ddots & \ddots & \vdots \\
	 &  & V_{N-1}^t & 1 
\end{array}\right),\quad
R_t(y)=
\left(\begin{array}{@{\,}cccc@{\,}}
	I_1^t & 1 &  &  \\
	 & I_2^t & \ddots &  \\
	 &  & \ddots & 1 \\
	y &  &  & I_{N}^t 
\end{array}\right),
\]
and $y$ is a complex variable. Let us introduce a new matrix
$X_t(y)$ defined by
\begin{equation}\label{ex}
X_t(y):=L_t(y)R_{t+M-1}(y)\cdots R_{t+1}(y)R_t(y).
\end{equation}
From (\ref{matrixform}) and (\ref{ex}), we obtain
\begin{equation}\label{Lax}
X_{t+1}(y)R_{t}(y)=R_t(y)X_t(y),
\end{equation}
which implies that the characteristic polynomial of $X_t(y)$ 
is invariant under the time evolution. 
Let $F(x,y):=\det{(X_t(y)-xE)}$ be the characteristic
polynomial of $X_t(y)$ ($E$ is the unit matrix).
Denote the affine curve defined by
$
F(x,y)=0
$
by $\widetilde{C}$, and its completion by $C$.
Of course,
$C$ is invariant as well under the time evolution.
This projective curve $C$ is called the
\textit{spectral curve} of the hpdToda.

\subsubsection{Properties of the spectral curve}

Now let us list the behaviour of $C$, following \cite{iwao} \S 2.
\begin{itemize}
\item on $C$, there exists a point $P:(x,y)=(\infty,\infty)$ around which
there exists a local coordinate $k$
such that $x=k^{-M}+\cdots$ and $y=k^{-N}+\cdots$.
\item on $C$, there exists a point $Q:(x,y)=(\infty,0)$
around which there exists a local coordinate $k$
such that $x=E k^{-1}+\cdots$ and $y= k^{N}+\cdots$,
where 
$E=(\prod_{n=1}^N{\prod}_{j=0}^{M-1}{I_n^j})\cdot\prod_{n=1}^N{V_n^0}$.
\item the $M$ points $A_j:(x,y)=(0,\, (-1)^{N}\prod_{n}{\!I_n^j})$,\ \ 
$j=0,1,\dots,M-1$\, lie on $C$.
\item the point $B:(x,y)=(0,\,\prod_{n}{\!\!V_n^t})$ lies on $C$.
\item The projection $p_x:C\ni(x,y)\mapsto x\in\PP^1$ is $(M+1):1$, and the projection
$p_y:C\ni(x,y)\mapsto y\in\PP^1$ is $N:1$.
\item $C$ has genus $g=\frac{(N-1)(M+1)-m+1}{2}$, where $m$ is the 
greatest common divisor of $N$ and $M$.
\end{itemize}
Hereafter we assume $C$ is smooth unless otherwise stated.

\subsubsection{The eigenvector mapping}

An \textit{isolevel set} $\tee_C$ is the set of matrices $X(y)$
(eq.(\ref{ex}))
associated with 
the spectral curve $C$. 
Now we construct a map from $\tee_C$ to $\mathrm{Pic}^{g+N-1}(C)$,
called the \textit{eigenvector mapping}, which plays a very important role
in the present method.

Let $X=X(y)$ be an element of $\tee_C$.
If $(x,y)\in \widetilde{C}$, there exists a complex $N$-vector 
$\vect{v}(x,y)$
such that $X(y)\vect{v}(x,y)=x\,\vect{v}(x,y)$, up to constant multiple.
Then there exists a Zariski open subset $C^\circ$ of $\widetilde{C}$
over which the morphism $C^\circ\ni(x,y)\mapsto\vect{v}(x,y)\in\PP^{N-1}$ is 
uniquely determined.
Moreover,
for a smooth $C$, this morphism can be extended uniquely over the whole $C$.
Denote this morphism by $\Psi_{X}:C\to\PP^{N-1}$.

The eigenvector mapping $\varphi_C:\tee_C\to\mathrm{Pic}^d(C)$\ $(d=g+N-1)$ is 
a map defined by the formula:
\[
\varphi_C(X)=\Psi_{X}^\ast(\mathcal{O}_{\PP^{N-1}}(1)),
\]
where $\mathcal{O}_{\PP^{N-1}}(1)$ is the invertible sheaf
of hyperplane sections over $\PP^{N-1}$.
Note that
it is nontrivial to prove $\varphi_C(X)\in\mathrm{Pic}^d(C)$
(see \cite{iwao} \S 2).

The role of the eigenvector mapping is to embed the set $\tee_C$
into $\mathrm{Pic}^d(C)$.
The following proposition is originally obtained in 
van Moerbeke, Mumford \cite{Mumford}.
\begin{prop}[\cite{Mumford}, thm.{}\,3]
The eigenvector mapping 
$
\varphi_C:\tee_C\to\mathrm{Pic}^d(C)
$
is an embedding.
\end{prop}
\bigskip

Although the definition of the eigenvector mapping is abstract,
we can
have an explicit formula to express $\varphi_C(X)$ in the present situation.
\begin{lemma}[\cite{iwao}, \S 2]\label{lem2.2}
Let $\vect{v}(x,y)=
\left(\begin{array}{@{\,}c@{\,}}
	g_1 \\
	\vdots \\
	g_N
\end{array}\right)$
be an eigenvector of $X(y)$ belonging to $x$
$(g_i=g_i(x,y),\ i=1,\dots,N)$. 
Then it follows that
$\varphi_C(X)=[(g_1/g_N)_\infty]$.
\end{lemma}

On the other hand, the divisor $(g_1/g_N)$ has the following expression 
(\cite{Mumford} prop.{}\,1):
\begin{equation}\label{eq2.5}
(g_1/g_N)=\mathcal{D}_1+(N-1)P-\mathcal{D}_2-(N-1)Q,
\end{equation}
where $\mathcal{D}_1$ and $\mathcal{D}_2$ are general and positive 
divisors of degree $g$.

Let $\mathfrak{d}(X):=\mathcal{D}_2$. Lemma \ref{lem2.2}
is rewritten as $\varphi_C(X)=[\mathfrak{d}(X)+(N-1)Q]$.

\subsection{Linearization theorem}\label{two-two}

Consider the $N\times N$ matrix $X_t(y)$ defined by (\ref{ex})
and the associated spectral curve $C$.
Let $\sigma$ and $\tau$ be the isomorphisms on $\tee_C$ defined by:
\begin{equation}\label{eq2.6}
\sigma(X_t(y))=SX_t(y)S^{-1},\quad \mu(X_t(y))=R_t(y)X_t(y)R_t(y)^{-1}=X_{t+1}(y),
\end{equation}
where
$S=
\left(\begin{array}{@{\,}cccc@{\,}}
	0 & 1 &  &  \\
	 & 0 & \ddots &  \\
	 &  & \ddots & 1 \\
	y &  &  & 0
\end{array}\right).
$
For the hpdToda equation (\ref{toda1}--\ref{toda2.5}, \ref{matrixform}),
$\sigma$ is the $n$-shift operator: $n\mapsto n+1$ and 
$\mu$ is the $t$-shift operator: $t\mapsto t+1$.

By calculating the divisors $\mathfrak{d}(\sigma(X_t))$ and 
$\mathfrak{d}(\mu (X_t))$,
we have the following theorem which illustrates the flow of the hpdToda equation
on $\mathrm{Pic}^d(C)$:

\begin{thm}[\cite{iwao}]\label{thm}
$(\mathrm{I})$:
Let $\mathcal{D}$ be the divisor
$
\mathcal{D}=P-Q.
$
Then the following 
diagram is commutative.
\[
\begin{array}{ccccc}
 &\tee_C& \to &\mbox{Pic}^d(C)& \\[2.5mm]
{\sigma}& \downarrow& &\downarrow&\hspace{-4mm}{+[\mathcal{D}]}
 \\[3mm]
  & \tee_C& \to& \mbox{Pic}^d(C)&
\end{array}.
\]
$(\mathrm{II})$:
Let $\mathcal{E}_{j}$ $(j=1,2,\dots,M)$ be the divisor
$
\mathcal{E}_j=P-A_j.
$ 
If $t\equiv j\pmod{M}$, the following 
diagram is commutative.
\[
\begin{array}{ccccc}
 &\tee_C& \to &\mbox{Pic}^d(C)& \\[2.5mm]
{\mu}& \downarrow& &\downarrow&\hspace{-4mm}{+[\mathcal{E}_j]} \\[2.5mm]
  & \tee_C& \to& \mbox{Pic}^d(C)&
\end{array}.
\]
\end{thm}
\begin{cor}
The time evolution $t\mapsto t+M$ is expressed as $Z\mapsto Z+[B-Q]$
on $\mathrm{Pic}^d(C)$.
\end{cor}
\proof
By theorem \ref{thm} (II), on $\mathrm{Pic}^d(C)$, 
$\{t\mapsto t+M\}$ is expressed by the formula: 
$Z \mapsto Z+[MP-A_0-A_1-\dots-A_{M-1}]$.
Then the relation $(x)=-MP-Q+A_0+A_1+\dots+A_{M-1}+B \in\mathrm{Div}^0(C)$ 
yields the result.
$\qed$

\begin{cor}\label{cor2.5}
The divisor $\mathcal{D}_1$ in $(\ref{eq2.5})$ satisfies $\mathcal{D}_1=\mathfrak{d}
(\sigma (X_t))$.
\end{cor}
\proof
By (\ref{eq2.5}),
$
[\mathcal{D}_1]=[\mathfrak{d}(X_t)+(N-1)Q-(N-1)P]
=[\mathfrak{d}(\sigma^{-N+1}(X_t))]
=[\mathfrak{d}(\sigma(X_t))]
$.
Because $\mathcal{D}_1$ and $\mathfrak{d}(\sigma(X_t))$ are general, positive 
and of degree $g$, it follows that $\mathcal{D}_1=\mathfrak{d}(\sigma(X_t))$.
$\qed$

\begin{cor}\label{cor2.6}
Let $\vect{v}(x,y)=
\left(\begin{array}{@{\,}c@{\,}}
	g_1 \\
	\vdots \\
	g_N
\end{array}\right)$
be an eigenvector of $X(y)$ belongs to $x$. 
Then
$(\mathrm{i})$ $(g_1/g_N)
=\mathfrak{d}(\sigma X)+(N-1)P-\mathfrak{d}(X)-(N-1)Q$, and\\
$(\mathrm{ii})$ 
$(g_{N}/yg_{N-1})=\mathfrak{d}(X)+(N-1)P-\mathfrak{d}(\sigma^{-1}X)-(N-1)Q$.
\end{cor}
\proof Part (i) follows immediately from (\ref{eq2.5}) and corollary \ref{cor2.5}.
Applying (i) to the matrix $\sigma^{-1}X=S^{-1}XS$ and noticing that 
$S\cdot (g_N\,y^{-1},g_1,\dots,g_{N-1})^T=(g_1,g_2,\dots,g_N)^T$,
we obtain (ii). $\qed$

\begin{rem}\label{rem2.1}
The time evolution $t\mapsto t+M$ is given by the map:
$\nu(X_t(y)):=L_t^{-1}(y)X_t(y)L_t(y)$.
In fact, $(\ref{ex},\ref{Lax})$ proves that $\nu(X_t(y))=X_{t+M}(y)$.
\end{rem}

\section{Tau function solution of the hpdToda equation}\label{sec3}

In this section, we assume $\mathrm{g.c.d.}(N,M)=1$.

\subsection{Construction of tau functions}

We construct a theta function solution of hpdToda equation.
As in the previous section, $X_t=X_t(y)$ denotes the
square matrix defined by (\ref{ex}).

Let $C$ be the (smooth) spectral curve associated with $X_t$.
Fix a symplectic basis $\alpha_1,\dots,\alpha_g;\beta_1,\dots,\beta_g$ of $C$
and the normalized holomorphic differentials $\omega_1,\dots,\omega_g$ such that
$\int_{\alpha_i}{\omega_j}=\delta_{i,j}$.
The $g\times g$ matrix $\Omega:=(\int_{\beta_i}{\omega_j})_{i,j}$ is called 
the \textit{period matrix} of $C$.
For a fixed point $p_0\in C$, the \textit{Abel-Jacobi mapping} 
$\abel:\mathrm{Div}(C)\to\CC^g/(\ZZ^g+\Omega\ZZ^g)$ is the
homomorphism defined by:
\[\textstyle
\sum{Y_i}-\sum{Z_j}
\ \mapsto\  
\sum(\int_{p_0}^{Y_i}{\omega_1},\cdots,
\int_{p_0}^{Y_i}{\omega_g})-\sum(\int_{p_0}^{Z_j}{\omega_1},\cdots,
\int_{p_0}^{Z_j}\omega_g).
\]

Let us consider the universal covering $\pi:\mathfrak{U}\to C$ and
fix an inclusion $\iota:C\hookrightarrow \mathfrak{U}$.
For simplicity, we slightly abuse the notation ``$\pi$" and
``$\iota$" to express the derived maps
$\mathrm{Div}(\mathfrak{U})\to\mathrm{Div}(C)$ and
$\mathrm{Div}(C)\hookrightarrow\mathrm{Div}(\mathfrak{U})$,
respectively.
Naturally, there exists a
continuous lift $\Abel:\mathrm{Div}(\mathfrak{U})\to \CC^g$
such that $\Abel\circ\iota(p_0)=0$.
For the projection $\rho:\CC^g\to\CC^g/(\ZZ^g+\Omega\ZZ^g)$,
it follows that $\rho\circ\Abel=\abel\circ\pi$.

\bigskip

For fixed $t\in\ZZ$, assume that some lifted positive divisor $\mathfrak{D}(X_t)\in
\mathrm{Div}^g(\mathfrak{U})$ with 
$\pi(\mathfrak{D}(X_t))=\mathfrak{d}(X_t)$ is specified.
Then there uniquely exist two positive divisors $\mathfrak{D}(\sigma X_t),\,
\mathfrak{D}(\mu X_{t})\in\mathrm{Div}^g(\mathfrak{U})$ such that:
\begin{eqnarray}
\Abel(\mathfrak{D}(\sigma X_t))=\Abel(\mathfrak{D}(X_t)+\iota P-\iota Q),\quad
\pi(\mathfrak{D}(\sigma X_t))=\mathfrak{d}(\sigma X_t),\\
\Abel(\mathfrak{D}(\mu X_{t}))=\Abel(\mathfrak{D}(X_t)+\iota P-\iota A_j),\quad
\pi(\mathfrak{D}(\mu X_{t}))=\mathfrak{d}(\mu X_{t}),
\end{eqnarray}
where $ t\equiv j\pmod{M}$.

Let $\tau^t$ be a
holomorphic function over $\mathfrak{U}$ defined by the formula:
\begin{equation}
\textstyle
\tau^t(p)=\theta
\left(
\Abel\{\mathfrak{D}(X_t)-p-\iota\Delta\}
\right),\qquad p\in \mathfrak{U},
\end{equation}
where $\theta(\bullet)=\theta(\bullet ;\Omega)$ is the Riemann theta function
and $\Delta\in\mathrm{div}^{g-1}(C)$ is the theta characteristic divisor of $C$
(\cite{tata}, Chap.{}\,II, cor.{}\,3.11).
To avoid cumbersome notations,
we often omit the letters ``$\Abel$", ``$\iota$" and use a simpler expression
$\tau^{t}(p)=\theta(\mathfrak{D}(X_t)-p-\Delta)$ 
when there is no confusion possible.

Although defined over $\mathfrak{U}$,
$\tau^t(p)$ can also be thought of as a multi-valued holomorphic function over 
$C$.
By the Riemann vanishing theorem (\cite{tata}, Chap.{}\,II, thm.{}\,3.11), 
the zero divisor of
$\tau^{t}(p)$ corresponds with $\mathfrak{d}(X_t)$.
\bigskip

Let $\tau^t_+(p):=\theta(\mathfrak{D}(\sigma X_t)-p-\Delta)$.
Then, by theorem \ref{thm}, the function
\[
\Psi^t(p):=
\frac{\tau_+^t(p)\cdot \tau^{t+1}(p)}
{\tau^t(p)\cdot \tau_+^{t+1}(p)}
=\frac{\theta(\mathfrak{D}(\sigma X_t)-p-\Delta)
\cdot\theta(\mathfrak{D}(\mu X_t)-p-\Delta)}
{\theta(\mathfrak{D}(X_t)-p-\Delta)\cdot
\theta(\mathfrak{D}(\mu\sigma X_t)-p-\Delta)}
\]
satisfies $[(\mbox{the zeros of denominator})]=
[(\mbox{the zeros of numerator})]\in\mathrm{Pic}^{2g}(C)$ and
therefore, it is a single-valued and meromorphic function over $C$.

Consider an eigenvector 
$X_t(y)
\left(\begin{array}{@{\,}c@{\,}}
	g_1^t \\
	\vdots \\
	g_N^t
\end{array}\right)
=
x
\left(\begin{array}{@{\,}c@{\,}}
	g_1^t \\
	\vdots \\
	g_N^t
\end{array}\right)
$,\ \  
$(g_i^t=g_i^t(x,y)=g_i^t(p))$.
From the relation $(g_1^t/g_N^t)
=\mathfrak{d}(\sigma X_t)+(N-1)P-\mathfrak{d}(X_t)-(N-1)Q$
(corollary \ref{cor2.6})
we derive the following equation by means of Liouville's theorem:

\begin{equation}\label{eq3.2}
\Psi^t(p)=c\times\frac{g_1^t(p)\cdot g_N^{t+1}(p)}{g_N^t(p)\cdot g_1^{t+1}(p)},
\qquad c:\mbox{constant}.
\end{equation}

By virtue of (\ref{eq3.2}), we can calculate some special values of $\Psi^t(p)$:
\begin{lemma}\label{lemma3.1}
On condition that $\mathrm{g.c.d}(N,M)=1$, we have
$(\mathrm{i})$
$\displaystyle
\Psi^t(P)=c$,
$(\mathrm{ii})$ $\displaystyle \Psi^t(Q)=c\times \frac{I_N^t}{I_1^t}$.
\end{lemma}
\proof
The lemma is proved by an elementary calculation, which we shall give in the appendix.
$\qed$

Because $\theta(\mathfrak{D}(X)-\iota Q-\Delta)=
\theta(\mathfrak{D}(X)+(\iota P-\iota Q)-\iota P-\Delta)
=\theta(\mathfrak{D}(\sigma X)-\iota P-\Delta)$,
it follows that
\[
\Psi^t(Q)=\Psi^t_+(P),\qquad\mbox{where}\quad \Psi^t_+(p)=
\frac{\tau^t_{++}(p)\cdot \tau^{t+1}_+(p)}{\tau^t_+(p)\cdot \tau_{++}^{t+1}(p)}.
\]
Then lemma \ref{lemma3.1} implies $I_1^t\Psi^t_+(P)=I_N^t\Psi^t(P)$.

Repeating this argument for $\Psi_+(p)$, we also derive
$I_2^t\Psi^t_{++}(P)\!=\!I_1^t\Psi_+^t(P)$, and inductively,
we have that:
\[
I_N^t\Psi^t(P)=I_1^t\Psi^t_+(P)=I_2^t\Psi^t_{++}(P)=I_3^t\Psi^t_{+++}(P)=\cdots.
\]
Let $\Psi_n^t:=\Psi^t_{++\cdots +}(P)$ ($n$ ``$+$"s).
Finally we obtain the equations 
$\Psi_{n+N}^t=\Psi_n^t$ and
$I_n^t\Psi_n^t=d$,
where the number $d$ does not depend on $n$.
\bigskip

Next consider the following single-valued meromorphic function over $C$:
\[
\Phi^t(p):=
\frac{\tau^t(p)\cdot \tau^{t+M}(p)}
{\tau^t_+(p)\cdot \tau_-^{t+M}(p)}
=\frac{\theta(\mathfrak{D}(X_t)-p-\Delta)
\cdot\theta(\mathfrak{D}(\nu X_t)-p-\Delta)}
{\theta(\mathfrak{D}(\sigma X_t)-p-\Delta)\cdot
\theta(\mathfrak{D}(\nu\sigma^{-1} X_t)-p-\Delta)}.
\]
Using corollary \ref{cor2.6} and Liouville's theorem,
we derive the following expression:
\begin{equation}\label{eq3.5}
\Phi^t(p)=c'\times\frac{g_N^t(p)\cdot g_N^{t+M}(p)}{g_1^t(p)\cdot g_{N-1}^{t+M}(p)
\cdot y},\qquad c':\mbox{constant},
\end{equation}
which again allows us to compute some special values of $\Phi^t(p)$.
\begin{lemma}\label{lemma3.2}
On condition that $\mathrm{g.c.d}(N,M)=1$, we have
$(\mathrm{i})$ $\displaystyle
\Phi^t(P)=c'$,
$(\mathrm{ii})$ $\displaystyle
\Phi^t(Q)=c'\times\frac{V_{N-1}^t}{V_N^t}$.
\end{lemma}
\proof See \ref{secA}. $\qed$

Due to $\Phi^t(Q)=\Phi^t_+(P)$ and lemma \ref{lemma3.2},
we have $V_N^t\Phi_+^t(P)=V_{N-1}^t\Phi^t(P)$, which implies
\[
V_{N-1}^t\Phi^t(P)=V_N^t\Phi_+^t(P)
=V_1^t\Phi_{++}^t(P)=V_2^t\Phi_{+++}^t(P)=\cdots.
\]
Let $\Phi_{n-1}^t:=\Phi_{++\cdots+}^t(P)$ ($n$ ``$+$"s).
Therefore we obtain 
$\Phi_{n+N}^t=\Phi_n^t$ and
$V_n^t\Phi_n^t=d'$,
where the number $d'$ does not depend on $n$.
\bigskip

Define 
$\tau_{-1}^t:=\tau^t(\iota P)$, 
$\tau_0^t:=\tau^t_+(\iota P)$, 
$\tau_1^t:=\tau_{++}^t(\iota P),
\cdots,
\tau_{n-1}^t:=\tau^t_{++\dots+}(\iota P)$ ($n$ ``$+$"s).
By the arguments above,
$I_n^t$ and $V_n^t$ have following expressions:
\begin{equation}\label{eq3.3}
I_n^t=d\times
\frac{\tau_{n-1}^t\cdot \tau_n^{t+1}}{\tau_n^t\cdot \tau_{n-1}^{t+1}},\qquad
V_n^t=d'\times\frac{\tau_{n+1}^t\cdot \tau_{n-1}^{t+M}}
{\tau_n^t\cdot\tau_n^{t+M}}.
\end{equation}

\subsection{Solution of hpdToda}

For $g$-dimensional vectors $\vect{a}$ and $\vect{b}$, 
$\langle \vect{a},\vect{b} \rangle$ denotes
$\vect{a}^T \vect{b}\in\CC$.

By periodicity $\mathfrak{d}(\sigma^N X_t)=\mathfrak{d}(X_t)$,
there exist integer vectors $\vect{n},\,\vect{m}\in\ZZ^g$ such that 
$
\Abel(N(\iota P-\iota Q))=\vect{n}+\Omega\vect{m}
$.
Considering the definition of the Riemann theta function 
(see \cite{tata}, \S II.1, for example),
we have 
\[
\tau^t_{n+N}=\tau_{n}^t\times\exp(-2\pi\mathrm{i}\cdot
\langle\vect{m},\vect{z}\rangle-\pi\mathrm{i}\cdot
\langle \vect{m},\Omega\vect{m}\rangle),\qquad
\mathrm{i}=\sqrt{-1},
\]
where $\vect{z}=\Abel(\mathfrak{D}(\sigma^{n+1}X_t)-\iota P-\Delta)$.
By (\ref{eq3.3}), we have
\begin{align}
&I_1^tI_2^t\cdots I_N^t=d^N\times\frac{\tau_1^t\cdot\tau_{N+1}^{t+1}}
{\tau_{N+1}^t\cdot\tau_1^{t+1}}
=d^N\times \exp(-2\pi\mathrm{i}\cdot
\langle  \vect{m}, \Abel(\iota P-\iota A_j) \rangle),\label{eq3.6}\\
&V_1^tV_2^t\cdots V_N^t={d'}^N\times
\frac{\tau_{N+1}^t\cdot\tau_{0}^{t+M}}{\tau_1^t\cdot\tau_{N}^{t+M}}\nonumber\\
&\ \ \ ={d'}^N\times\exp(-2\mathrm{i}\,\pi\cdot\langle
\vect{m},\Abel(\iota A_0+\cdots+\iota A_{M-1}-(M-1)\iota P-\iota Q)
\rangle),
\label{eq3.7}
\end{align}
where $j\equiv t\pmod{M}$.
Recall $\prod_n{I_n^{t+M}}=\prod_n{I_n^t}$ and $\prod_n{V_n^{t+1}}=\prod_n{V_n^t}$,
which imply that 
$d$ depends on $t\pmod{M}$ and that $d'$ is independent from $t$.
Finally we obtain the conclusion:

\begin{thm}
If $\mathrm{g.c.d.}(N,M)=1$,
$(\ref{eq3.3}\mbox{--}\ref{eq3.7})$ solves
the hpdToda $(\ref{toda1}\mbox{--}\ref{toda2.5})$.
\end{thm}

\section{The general cases}\label{sec4}
In the previous sections, we have assumed that $\mathrm{g.c.d.}(N,M)=1$.
Unfortunately, the method which we have established in this paper cannot be applied
in the general cases.

For example, when $N=M=2$, the characteristic polynomial of the matrix 
$X_t(y)$ (equation (\ref{ex}))
is: 
\[
\det{(X_t(y)-xE)}= y^2-y(2x+U_1)+x^2-U_2x+U_3-U_4y^{-1},
\]
where $U_1=I_1^tI_2^t+I_1^{t+1}I_2^{t+1}+V_1^tV_2^t$, $U_2=I_1^tI_1^{t+1}+I_2^tI_2^{t+1}
+I_1^tV_2^t+I_1^{t+1}V_1^t+I_2^tV_1^t+I_2^{t+1}V_2^t$,
$U_3=I_1^tI_2^tI_1^{t+1}I_2^{t+1}+I_1^{t+1}I_2^{t+1}V_1^tV_2^t+V_1^tV_2^tI_1^tI_2^t$,
$U_4=I_1^tI_2^tI_1^{t+1}I_2^{t+1}V_1^tV_2^t$.
However, the hungry Toda system (\ref{toda1}--\ref{toda2.5})
has the extra conserved quantity $I_1^t+I_2^t+I_1^{t+1}+I_2^{t+1}+V_1^t+V_2^t$,
which is independent from $U_1$, $U_2$, $U_3$ and $U_4$.
This means that the spectral curve does not faithfully reflect the data
of the system.

For this reason, we should try to trace the problem to the case $\mathrm{g.c.d.}(N,M)=1$.
Denote by $\mathrm{Toda}_{N,M}$ the hungry Toda system (\ref{toda1}--\ref{toda2.5})
associated with the positive integers $N$ and $M$.
It is sufficient to prove the following statement.

\begin{prop}\label{prop4.1}
Define the initial values $I_n^0:=\zeta+o(\zeta)$, $(\zeta\to\infty,\forall n)$
for some complex parameter $\zeta$,
and let $\{I_n^t,V_n^t\}_{n,t}$ be a solution of $\mathrm{Toda}_{N,M}$.
When $\zeta\to \infty$,
the new sequence 
\[
\{I_n^{kM+1},I_n^{kM+2},\dots,I_n^{kM+M-1},V_n^{kM+1},V_n^{kM+2},\dots,V_n^{kM+M-1}\}_{n,k}
\]
is a solution of $\mathrm{Toda}_{N,M-1}$.
\end{prop}
\proof
We shall prove the following:
\begin{align}
&I_n^{kM+M-1}=I_n^{kM-1}+V_n^{kM-1}-V_{n-1}^{kM+1}+o(1),\label{eq.17}\\
&V_{n}^{kM+1}=\frac{I_{n+1}^{kM-1}V_n^{kM-1}}{I_n^{kM+M-1}}\cdot
(1+o(1)).\label{eq.18}
\end{align}
By (\ref{toda1}--\ref{toda2.5}) and \textit{Remark} (page \pageref{rema}),
we have
\[
I_n^{t}=\zeta+o(\zeta),\ (\forall n)\quad \Rightarrow\quad 
\left\{
\begin{array}{l}
I_n^{t+M}=\zeta+o(\zeta),\ (\forall n)\\
V_n^{t+1}=V_n^{t}+o(1),\ (\forall n)
\end{array}
\right.\qquad (\zeta\to\infty).
\]
Then, in our situation, it follows that $V_n^{kM+1}=V_n^{kM}+o(1)$
for all $k\in\ZZ_{\geq 0}$ and $n$.
Using (\ref{toda1}--\ref{toda2.5}) again, we derive
equations (\ref{eq.17},\ref{eq.18}). $\qed$

Applying proposition \ref{prop4.1} repeatedly, 
we can trace the problem to the case $\mathrm{g.c.d.}(N,M)=1$.

\begin{example}
The hungry Toda system with $N=M=2$ can be traced to the case $N=2,M=3$.

Let 
$
L_0:=
\left(\begin{array}{@{\,}cc@{\,}}
	1 & V_2^0\,y^{-1} \\
	V_1^0 & 1
\end{array}\right)
$,
$
R_0:=
\left(\begin{array}{@{\,}cc@{\,}}
	\zeta & 1\\
	y & \zeta
\end{array}\right)
$,
$
R_1:=
\left(\begin{array}{@{\,}cc@{\,}}
	I_1^0 & 1 \\
	y & I_2^0
\end{array}\right)
$,
$
R_2:=
\left(\begin{array}{@{\,}cc@{\,}}
	I_1^1 & 1 \\
	y & I_2^1
\end{array}\right)
$.
Define $X_0:=L_0R_2R_1R_0$.
The characteristic polynomial of $X_0$ is:
\begin{eqnarray*}
\det{(X_0-xE)}=&-y^3+y^2(\zeta^2+U_1)-y\{(2\zeta+U_5)x+U_1\zeta^2+U_3\}\\
&+x^2-(U_2\zeta+U_6)x+U_3\zeta^2+U_4-U_4\zeta^2y^{-1},
\end{eqnarray*}
where $U_5=I_1^0+I_2^0+I_1^1+I_2^1+V_1^0+V_2^0$ and 
$U_6=I_1^0I_1^1V_1^0+I_2^0I_2^1V_2^0$.
Note that $U_5$ is the hidden conserved quantity of $\mathrm{Toda}_{2,2}$.
Let $\{I_n^t,V_n^t\}_{n,t}$ be the solution of $\mathrm{Toda}_{2,3}$.
Then the sequence 
\begin{eqnarray*}
\lim_{\zeta\to\infty}I_n^0,\,
\lim_{\zeta\to\infty}I_n^1,\,
\lim_{\zeta\to\infty}I_n^3,\,
\lim_{\zeta\to\infty}I_n^4,\,
\lim_{\zeta\to\infty}I_n^6,\,\dots;\\
\lim_{\zeta\to\infty}V_n^0,\,
\lim_{\zeta\to\infty}V_n^1,\,
\lim_{\zeta\to\infty}V_n^3,\,
\lim_{\zeta\to\infty}V_n^4,\,
\lim_{\zeta\to\infty}V_n^6,\dots.
\end{eqnarray*}
solves $\mathrm{Toda}_{2,2}$.
\end{example}

\subsection*{Acknowledgement}
The author is very grateful to Professor Tetsuji Tokihiro and
Professor Ralph Willox for helpful comments 
on this paper. 
This work was supported by KAKENHI 09J07090.

\appendix
\section{Proofs of lemmas}\label{secA}
Let $\Psi^t(p)$ and $\Phi^t(p)$ be the meromorphic functions
defined in section
\ref{sec3}.
We shall now prove lemma \ref{lemma3.1}, \ref{lemma3.2}.
In the appendix, we assume $\mathrm{g.c.d.}(N,M)=1$.

Denote the set of $N\times N$ matrices by $M_N(\CC)$ and
the subset of diagonal matrices by $\Gamma\subset M_N(\CC)$.
For a matrix $X\in M_N(\CC)$ and subsets $A,\,B\subset M_N(\CC)$,  
let 
$A+X:=\{a+X\,\vert\,a\in A\}$,
$A X:=\{a X\,\vert\,
a\in A\}$, 
$A+B:=\{a+b\,\vert\,a\in A,b\in B\}$ and 
$AB:=\{ab\,\vert\,a\in A,b\in B\}$.

For two meromorphic functions $f,g$ over $C$ and a point $p\in C$, 
``$f\sim g$ around $p$" means
$0<\lim_{z\to p}{\zet{f(z)/g(z)}}<+\infty$.
\bigskip

Let $(g_1,g_2,\dots,g_N)^T$ be an eigenvector of $X=X(y)\in\tee_C$
belonging to an eigenvalue $x$. Then
$g_1,\dots,g_{N}$ are meromorphic functions over $C$.
The following lemma is fundamental. 
\begin{lemma}\label{lemmaa.1}
$(\mathrm{i})$ Let $k$ be a local coordinate around $P$.
Then $g_1/g_N= k^{N-1}+\cdots$, 
$g_2/g_N= k^{N-2}+\cdots$, $\dots$, $g_{N-1}/g_N=k+\cdots$.\\
$(\mathrm{ii})$  Let $k$ be a local coordinate around $Q$.
Then $g_1/g_N\sim k^{-N+1}$, $g_2/g_N\sim k^{-N+2}$, $\dots$, 
$g_{N-1}/g_N\sim k^{-1}$.
\end{lemma}
\proof (i) 
Recall that we have $x=k^{-M}+\cdots$ and $y=k^{-N}+\cdots$ around $P$. 
By (\ref{ex}), $X_t$ is contained in the subset
$(E+\Gamma S^{-1})(\Gamma+S)^M=\Gamma S^{-1}+\Gamma+\Gamma S
+\dots+\Gamma S^{M-1}+S^{M}$.
Then the equation $X_t(y)\,\vect{v}=x\,\vect{v}$ implies:
\[
(\gamma_{-1}S^{-1}+\gamma_0+\gamma_1S+\dots+\gamma_{M-1}S^{M-1}+S^M)\cdot\vect{v}
=k^{-M}\vect{v}+(\mbox{higher terms}),
\]
where $\gamma_i$ $(i=-1,0,\dots,M-1)$ are diagonal matrices.
Let $T:=kS$.
Therefore we obtain
$\left(T^M+\sum_{i=-1}^{M-1}{k^{M-i}\gamma_i T^i}\right)
\cdot\vect{v}=\vect{v}+(\mbox{higher})$.
Because $N$ and $M$ are relatively prime,
the solution of $T\vect{v}=\vect{v}$ is
$\vect{v}=(k^{N-1},k^{N-2},\dots,1)^{T}$ up to a constant multiple.
This fact leads to the desired result.\\
(ii)
Let $k$ be a local coordinate around $Q$ such that $x=Ek^{-1}+\cdots$ and 
$y=k^{M}+\cdots$ (Section \ref{sec2}). 
It follows that
\[
(\gamma_{-1}S^{-1}+\gamma_0+\gamma_1S+\dots+\gamma_{M-1}S^{M-1}+S^M)\cdot\vect{v}
=Ek^{-1}\vect{v}+(\mbox{higher}).
\]
Let $U:=k^{-1}S$.
Then we have
$\left(\gamma_{-1}U^{-1}+\sum_{i=0}^{M}{k^{i+1}\gamma_i U^i}\right)
\cdot\vect{v}=E\vect{v}+(\mbox{higher})$.
Standard results from linear algebra prove that
there exist $(N-1)$ complex numbers $c_1,\dots,c_{N-1}$ such that
\[
U\cdot(c_1k^{-N+1},c_2k^{-N+2},\dots,1)^{T}=E\cdot(c_1k^{-N+1},c_2k^{-N+2},\dots,1)^T,
\]
which leads to the desired result.$\qed$

\subsection*{Proof of lemma \ref{lemma3.1}}

The equation $X_{t+1}(y)R_t(y)=R_t(y)X_t(y)$ (\ref{matrixform}) implies
$(g^{t+1}_1,g_2^{t+1}\dots,g^{t+1}_N)=R_t(y)\cdot
(g^{t}_1,g_2^{t}\dots,g^{t}_N)$.
Then (\ref{eq3.2}) gives rise to
\[
\Psi^t(p)=c\times
\frac{g_1^t}{g_N^t}\cdot\frac{I_N^tg^t_N+g^t_1y}{I_1^tg^t_1+g^t_2}.
\]
By lemma \ref{lemmaa.1}, $\Psi^t$ satisfies $\Psi^t=c+\cdots$, around $P$,
and $\Psi^t=c\cdot(I_N^t/I_1^t)+\cdots$, around $Q$.
$\qed$

\subsection*{Proof of lemma \ref{lemma3.2}}

As mentioned in remark \ref{rem2.1}, one has that 
$L_t(y)X_{t+M}(y)=X_t(y)L_t(y)$, which implies
$(g^{t}_1,g_2^{t}\dots,g^{t}_N)=L_t(y)\cdot
(g^{t+M}_1,g_2^{t+M}\dots,g^{t+M}_N)$.
Then (\ref{eq3.5}) leads
\[
\Phi^t(p)=c'\times\frac{V_{N-1}^{t}\,g_{N-1}^{t+M}+g_N^{t+M}}
{V_N^t\,g_N^{t+M}y^{-1}+g_1^{t+M}}\cdot
\frac{g_N^{t+M}}{g_{N-1}^{t+M}\cdot y}.
\]
By lemma \ref{lemmaa.1}, $\Phi^t$ satisfies $\Phi^t=c'+\cdots$, around $P$,
and $\Phi^t=c'\cdot(V_{N-1}^t/V_N^t)+\cdots$, around $Q$.
$\qed$

\end{document}